\input amstex
\documentstyle{amsppt}
\magnification=\magstep1
\NoRunningHeads
\topmatter
\title
 Topological Rigidity for Certain Families of Differentiable Plane Curves
\endtitle
\author Jean Paul DUFOUR and Yasuhiro KUROKAWA
\endauthor

\subjclass 
58C27
\endsubjclass

\abstract
We show that the topological classification and the smooth classification
are generically the same for certain families of plane curves
in a semi-local case(the double local case).
Especially we give the normal form of transversely jointed two families 
of plane curves with second order contact at the envelope.

\medpagebreak

\noindent
{\it Keywords} : singularity, divergent diagram, semi-local,
web structure, rigidity

\endabstract
\endtopmatter
\document
\baselineskip=19pt
\head 1. Introduction
\endhead
The divergent diagrams of map germs
$(\Bbb R, 0)\leftarrow (\Bbb R ^{2}, 0)\rightarrow(\Bbb R ^{2}, 0)$
appear in several 
subjects in geometry(e.g.[1, 2, 4, 6, 7, 8, 10, 13]).
For these subjects, in any case such divergent 
diagrams are identified up to smooth coordinate changes in each spaces.
Particularly in [1, 4] a generic classification of such divergent diagrams
have been studied.
Moreover in [4] it was shown that their topological classification
and their $C^{\infty}$-classification are generically the same.
Such property is called {\it topological rigidity}.
This classification have been applied to the vision theory([6, 9, 10]).

In a semi-local case of the above studies 
the following diagram of map germs appear:
$$\CD
(\Bbb R,0) @<<< (\Bbb R^{2},0) @.    @.\\
\;\;\;\;\;\; @. @. @. @. @. @.\;\;\;\;\;\;\;\;\;\;\;\;\;\; (\Bbb R^2, 0)\\
(\Bbb R,0) @<<< (\Bbb R^{2},0) @.    @.  \\
\endCD$$
Generic types of such diagrams were studied in [6, 10] and certain
normal forms were given in [10] by using the standard method of
singularity theory. 
However their topological rigidity is unsolved problem.
In this paper we discuss the following problem.

The two diagrams of $C^{\infty}$-map germs
$(f_1^{j}, \gamma_1^{j}; f_2^{j}, \gamma_2^{j}), j=1,2$ :
$$\CD
(\Bbb R,0) @<f_1^{j}<< (\Bbb R^{2},0) @>\gamma_1^{j} >> (\Bbb R^2,0) 
@<\gamma_2^{j} << (\Bbb R^{2},0) @>f_2^{j}>> (\Bbb R,0) 
\endCD$$
are called $C^{\infty}$-(resp.{\it topologically\/}){\it equivalent}
if there exists $C^{\infty}$-diffeomorphism germs (resp. homeomorphism
germs)
$h_i : (\Bbb R , 0)\to (\Bbb R , 0), 
H_i : (\Bbb R^2 , 0)\to (\Bbb R^2 , 0), i=1, 2$ and
$K : (\Bbb R^2 , 0)\to (\Bbb R^2 , 0)$ such that
$h_i\circ f_i^1 = f_i^2\circ H_i,\, K\circ\gamma_i^1
= \gamma_i^2 \circ H_i$ for $i=1,2.$

\noindent
{\it Topological Rigidity Problem} :
For the diagrams of map germs
$$(\Bbb R, 0)\leftarrow (\Bbb R ^{2}, 0)\rightarrow(\Bbb R ^{2}, 0)
\leftarrow (\Bbb R^2 , 0)\rightarrow (\Bbb R , 0)$$
the topological classification and the $C^{\infty}$-classification
are generically the same or not? 

Our result states that the answer of this problem is affirmative.
More precisely to state our theorem we shall recall some results in the following section.

\head 2. Previous results and statements of Theorems
\endhead
At first we shall recall some results of families of plane curves in [1, 3, 4, 6, 10, 14].
Denote by $\Cal E _{x_{1},\ldots , x_{n}}$ (or briefly $\Cal E _n$)
the ring of all smooth function germs
on $\Bbb R ^{n}$ at $0$ with coordinates $(x_{1},\ldots , x_{n})$
and denote by $\Cal M _{x_{1},\ldots , x_{n}}$ (or briefly $\Cal M_n$)
the unique maximal ideal of
$\Cal E _{x_{1},\ldots , x_{n}}.$
Denote by $S_{f}$ the singular set of a $C^{\infty}$-map germ $f:(\Bbb R^n, 0)\to \Bbb R^m.$
We shall suppose that all map germs are of class $C^{\infty}$
unless otherwise stated.

In [4] the generic type of $(\Bbb R, 0)\leftarrow (\Bbb R ^{2}, 0) \rightarrow(\Bbb R ^{2}, 0)$ have been given as follows:

\noindent$(I)$ $f$ is a submersion and $\gamma$ is regular.

\noindent$(II)$ $f$ is of Morse type and $\gamma$ is regular.

\noindent$(III)$ $f$ is a submersion,
$\gamma$ is a fold,
$f$ restricted to the singular set $S_\gamma$ of $\gamma$ is regular

and $(f, \gamma):(\Bbb R^{2},0)\to (\Bbb R^{3},0)$ is regular.  

\noindent$(IV)$ $f$ is a submersion,
$\gamma$ is a fold,
$f\vert_{S_{\gamma}}$  
is of Morse type

and $(f, \gamma):(\Bbb R^{2},0)\to (\Bbb R^{3},0)$ is regular.  

\noindent$(V)$ $f$ is a submersion, $\gamma$ is a fold,
$(f, \gamma):(\Bbb R^{2},0)\to (\Bbb R^{3},0)$ is a Whitney umbrella whose 

line of double points is transversal at $0$
to the direction $\{0\}\times\Bbb R^{2}$ in $\Bbb R^{3}$. 

\noindent$(VI)$ $f$ is a submersion,
$\gamma$ is a cusp
and $(f, \gamma):(\Bbb R^{2},0)\to (\Bbb R^{3},0)$ is regular. 

\medpagebreak

In the semi-local case generic types have been given as follows([6, 10]):

\noindent$(I, I)^{0}, (I, I)^{1}, (I, I)^{2}$ : 
$(f_{1}, \gamma_{1}), (f_{2}, \gamma_{2})$ are both of type (I)

and $(f_{1}\circ {\gamma_{1}}^{-1}, f_{2}\circ {\gamma_{2}}^{-1}): (\Bbb R ^{2}, 0)
\to (\Bbb R ^{2}, 0)$ is respectively regular, a fold, a cusp.

\medpagebreak
\noindent$(II, I)$ : 
$(f_{1}, \gamma_{1})$ is of type (II),
$(f_{2}, \gamma_{2})$ is of type (I)

and $(f_{1}\circ {\gamma_{1}}^{-1}, f_{2}\circ {\gamma_{2}}^{-1}) : (\Bbb R ^{2}, 0)
\to (\Bbb R ^{2}, 0)$ is a fold.

\medpagebreak
\noindent$(III, I)^{0} , (III, I)^{1}$ : 
$(f_{1}, \gamma_{1})$ is of type (III),
$(f_{2}, \gamma_{2})$ is of type (I)
and respectively 

$\gamma_1(S_{\gamma_{1}}) \, \bar\pitchfork \, \gamma_{2}(f_{2}^{-1}(0))$,
$\gamma_1(S_{\gamma_{1}}) \, \bar\pitchfork \, \gamma_{2}(f_{2}^{-1}(0))$ 
with two point contact.





\medpagebreak

\noindent$(IV, I)$ : 
$(f_{1}, \gamma_{1})$ is of type (IV),
$(f_{2}, \gamma_{2})$ is of type (I),
$\gamma_1(S_{\gamma_{1}}) \, \bar\pitchfork \, \gamma_{2}(f_{2}^{-1}(0))$ 

and $\bigcup\Sb{t\in (\Bbb R, 0)}\endSb\gamma_1 (S_{\lambda_{t}})\,\bar\pitchfork\,
\gamma_2 (f_2 ^{-1}(0))$
where $\lambda_{t}\in\Cal E _1$ such that $graph(\lambda_t )=f_1 ^{-1}(t).$

\medpagebreak
\noindent$(V, I)$ : 
$(f_{1}, \gamma_{1})$ is of type (V),
$(f_{2}, \gamma_{2})$ is of type (I),
$\gamma_1(S_{\gamma_{1}}) \, \bar\pitchfork \, \gamma_{2}(f_{2}^{-1}(0)),$

 $\bigcup\Sb{t\in (\Bbb R, 0)}\endSb\gamma_1 (S_{\lambda_{t}})\,\bar\pitchfork\,
\gamma_2 (f_2 ^{-1}(0))$
where $\lambda_{t}\in\Cal E _1$ such that $graph(\lambda_t )=f_1 ^{-1}(t)$

and the tangent cone of 
$\gamma_1(f_1 ^{-1}(0)) \, \bar\pitchfork \, \gamma_{2}(f_{2}^{-1}(0)).$ 

\medpagebreak
\noindent$(VI, I)$ : 
$(f_{1}, \gamma_{1})$ is of type (VI),
$(f_{2}, \gamma_{2})$ is of type (I)

and the tangent cone of $\gamma_1(S_{\gamma_{1}}) \; \bar{\pitchfork}
\; \gamma_{2}(f_{2}^{-1}(0)).$

\medpagebreak
\noindent$(III, III)$ :  
$(f_{1}, \gamma_{1}), (f_{2}, \gamma_{2})$ are both of type (III)
and $\gamma_1(S_{\gamma_{1}}) \, \bar\pitchfork \, \gamma_2(S_{\gamma_{2}}).$

\definition{Remark 1.1}
In case $(IV, I)$ a condition $\bigcup\Sb{t\in (\Bbb R, 0)}\endSb\gamma_1 (S_{\lambda_{t}})\,\bar\pitchfork\,
\gamma_2 (f_2 ^{-1}(0))$
is added to the previous generic condition of $(IV, I)$ in [10].
Even if this condition is added, the above condition of $(IV, I)$ is generic.
Similarly in case $(V, I)$ two conditions such that
$\bigcup\Sb{t\in (\Bbb R, 0)}\endSb\gamma_1 (S_{\lambda_{t}})\,\bar\pitchfork
\allowmathbreak
\gamma_2 (f_2 ^{-1}(0))$
and
the tangent cone of 
$\gamma_1(f_1 ^{-1}(0)) \, \bar\pitchfork \, \gamma_{2}(f_{2}^{-1}(0))$ 
are added in this paper.
These added conditions are effective to solve the rigidity problem. 
\enddefinition

The following figures  give families of plane curves 
$\gamma_i (f_i ^{-1}(t)), i=1,2$ for each generic types.

\newpage

.

\newpage

\topspace{9.3cm}

Normal forms of the generic types have been given in [1, 3, 4, 10, 14](See sect.5).
However we give another normal form for type (III, III)
which is more precise than the previous normal form([10])
and moreover effective to show the rigidity problem.

\proclaim{Theorem A}
The normal form for type (III, III) is the following:
$$f_{1}(x_1, y_1)=x_{1}+y_{1}, \gamma_{1}=(x_{1}, y_{1}^{2})$$
$$f_{2}(x_2, y_2)=x_{2}+y_{2}+\theta (x_{2}, y_{2}), \gamma_{2}=(x_{2}^{2}, y_{2})$$

\noindent
where $\theta$ is an arbitrary function germ in  $\Cal M _{x_{2}, y_{2}}$ satisfying $\theta (x_{2}, 0)=\theta (0, y_{2})=0.$
\endproclaim

\definition{Remark 1.2}
In [10] normal form of (III, III) was given  as follows: 
$f_{1}=y_{1}+\alpha _{1}\circ \gamma_{1},\;
\gamma_{1}=(x_1, y_{1}^{2});\;
f_{2}=x_{2}+\alpha _{2}\circ \gamma_{2},\;
\gamma_{2}=(x_{2}^{2}, y_{2}),$
where $\alpha _{1}, \alpha _{2}\in\Cal M _{u, v}$
with
$\frac{\partial\alpha_{1}}{\partial u}(0)\not=0,
\frac{\partial\alpha_{2}}{\partial v}(0)\not=0.$
\enddefinition

\definition{Remark 1.3}
In the semi-local case we remark that the followings: 
(1) The generic types are not topologically equivalent one another.
(2) The normal forms depend on arbitrary functions with some conditions,
that is so-called ``functional moduli'' appear in the normal forms
except (I, I)${}^{0}$, (I, I)${}^{1}$, (I, I)${}^2$ and (II, I).
(On the normal forms of each generic types, see sect. 5)
\enddefinition

Therefore to solve the topological rigidity problem it is sufficient to
consider each generic types with functional moduli.

\proclaim{Theorem B}
The topological rigidity holds for every generic types.
\endproclaim

\head{3. Proof of Theorem A}
\endhead
Firstly we detect the normal form for $\gamma_1$ and $\gamma_2.$
Next we shall find out certain coordinate transformations preserving
$\gamma_1$ and $\gamma_2$ which give the normal form in Theorem A.

\proclaim{Lemma 3.1}
Let
$\CD
(\Bbb R^{2}, 0) @>\gamma_{1}>> (\Bbb R ^{2}, 0) @<\gamma_{2}<< (\Bbb R^{2}, 0)
\endCD$
be a bi-germ such that $\gamma_{1}, \gamma_{2}$ are both fold map germs
and that the discriminant sets of $\gamma_{1}, \gamma_{2}$ are transversal at $0$ each other.
Then we can express the bi-germ for some coordinates as follows:
$$\gamma_{1}(x_{1}, y_{1})=(x_{1}, y_{1}^{2}),
\,  \gamma_{2}(x_{2}, y_{2})=(x_{2}^{2}, y_{2}).$$ 
\endproclaim
\demo{Proof}
We choose coordinates $(x_{1}, y_{1})$ in the source space of $\gamma_{1}$,
$(u, v)$ in the target space of $\gamma_{1}$ and $\gamma_{2}$ such that
$\gamma_{1}(x_{1}, y_{1})=(x_{1}, y_{1}^{2}).$
Then we can take coordinates $(x_2, y_2)$ such that 
$\gamma_{2}(x_{2}, y_{2})=(u\circ \gamma_{2}, y_{2}).$
From the Malgrange preparation theorem,
there exist $A, B \in \Cal E _{u, v}$ such that
$x_{2}^{2}=A\circ \gamma_{2}+2(B\circ \gamma_{2})x_{2}.$
Consider coordinate transformations 
$\phi_{1}, \psi , \phi_{2}$ 
which preserve $\gamma_1$
by
$\phi_{1}(x_{1}, y_{1})
=(A\circ \gamma_{1}(x_{1}, y_{1})+B^{2}\circ \gamma_{1}(x_{1}, y_{1}), y_{1}),$
$$\psi(u, v)
=(A(u, v)+B(u, v)^{2}, v),$$
$$\phi_{2}(x_{2}, y_{2})
=(x_{2}-B(u\circ \gamma_{2}(x_{2}, y_{2}), y_{2}), y_{2}).$$
By these coordinate changes we obtain the normal form of the bi-germ.
\qed
\enddemo

From now on we suppose that in this section always the bi-germ
$$\CD
(\Bbb R^{2}, 0) @>\gamma_{1}>> (\Bbb R ^{2}, 0) @<\gamma_{2}<< (\Bbb R^{2}, 0)
\endCD$$
is defined by
$\gamma_{1}(x_{1}, y_{1})=(x_{1}, y_{1}^{2}),
\,  \gamma_{2}(x_{2}, y_{2})=(x_{2}^{2}, y_{2}).$

The following Lemma is obvious:
\proclaim{Lemma 3.2}
Let
$\CD
(\Bbb R^{2}, 0) @>\gamma_{1}>> (\Bbb R ^{2}, 0) @<\gamma_{2}<< (\Bbb R^{2}, 0)
\endCD$
be a bi-germ given by the above.
Let
$H_{1}, H_{2}, K:(\Bbb R^{2},0)\to(\Bbb R^{2},0)$
be diffeomorphism germs defined by
$$H_{1}(x_{1}, y_{1})=(x_{1}A^{2}\circ \gamma_{1}, y_{1}B\circ \gamma_{1}),$$
$$H_{2}(x_{2}, y_{2})=(x_{2}A\circ \gamma_{2}, y_{2}B^{2}\circ \gamma_{2}),$$
$$K(u, v)=(uA^{2}, vB^{2}),$$
where $A, B\in \Cal E _{u, v}$
with $A(0)\not=0, B(0)\not=0$.
Then 
$\gamma_1\circ H_{1}=K\circ \gamma_1$ and $\gamma_2\circ H_{2}=K\circ \gamma_2$
hold.
\endproclaim

We say that the triple of the above diffeomorphism germ $(H_{1}, K, H_{2})$
is {\it a} $(\gamma_{1}, \gamma_{2})$-{\it compatible diffeomorphism germ}.
\proclaim{Lemma 3.3}
The germs of type (III, III) are equivalent to
$$(\Bbb R, 0) @<x+y<< (\Bbb R ^{2}, 0) @>\gamma_{1}>> (\Bbb R ^{2}, 0)
              @<\gamma_{2}<< (\Bbb R^{2}, 0) @>f>> (\Bbb R, 0),$$
where $f (x, y)=x+y+$higher term.
\endproclaim
\demo{Proof} Let $(f_{1}, \gamma_1 ; f_{2}, \gamma_{2})$ be type (III, III).
Firstly we shall show that the type (III, III) is equivalent to
$(x_1+y_1, \gamma_{1}; f, \gamma_2),$ where
$f \in \Cal M_{x, y}$ with $\frac{\partial f}{\partial x}(0)
\frac{\partial f}{\partial y}(0)\not= 0.$
Since $(f_{1}, \gamma_{1})$ is of type (III), 
$\frac{\partial f_{1}}{\partial y_{1}}(0)\not= 0.$
So by using the inverse function of $f_1 (0, y_1 )$ as a coordinate
change in the target of $f_1 ,$ we can suppose that
$f_{1}(0, y_{1})=y_{1}.$
Then by the Malgrange preparation theorem,
there exist $\alpha, \beta \in \Cal E _{2}$ such that
$f_{1}(x_{1}, y_{1})=\alpha(x_{1}, y_{1}^{2})+y_{1}\beta(x_{1}, y_{1}^{2}).$
Thus $\alpha(0, y_{1}^{2})=0.$
Moreover since $f_{1}\vert_{S_{\gamma_{1}}}$ is non-singular we may suppose
that $\frac{\partial\alpha}{\partial u}(0)>0.$
Hence we have $\alpha(x_{1}, y_{1}^{2})=x_{1}\widetilde{\alpha}
(x_{1}, y_{1}^{2})$ for some $\widetilde{\alpha}\in\Cal E_2$
with $\widetilde{\alpha}(0)>0.$
Now set $A=\root\of{\widetilde{\alpha}}, B=\beta$.
Then using the $(\gamma_{1}, \gamma_{2})$-compatible diffeomorphism 
given in Lemma 3.2
we obtain the required form, where
$f=f_{2}\circ H_2^{-1}.$

Next put $a={(\frac{\partial f}{\partial x}(0)/
        \frac{\partial f}{\partial y}(0))}^{-\frac{2}{3}},
b={(\frac{\partial f}{\partial y}(0))}^{-1}
{(\frac{\partial f}{\partial x}(0)/
        \frac{\partial f}{\partial y}(0))}^{-\frac{4}{3}}.$
Then consider the coordinate transformations
$h_{1}(t_{1})=at_{1}$ in the target of $f_{1}=x_{1}+y_{1}$,
$k(t)=bt$ in the target of $f$,
$(\gamma_{1}, \gamma_{2})$-compatible diffeomorphisms
$H_{1}, K, H_{2}$ by putting
$A=\root\of{a}, B=a$ in Lemma 3.2.
Using these coordinate transformations 
we obtain the required form in this Lemma.
\qed
\enddemo

In order to find coordinate transformations which reduce to
the normal form in Theorem A, we now need the following observation.
Suppose that the diagram
$$\CD
(\Bbb R,0) @<x+y << (\Bbb R^{2},0) @>\gamma_{1}>>
(\Bbb R^2,0) @<\gamma_{2}<< (\Bbb R^{2},0) @>f>>
(\Bbb R, 0)\\
 @V h VV         @V H_{1} VV        @V K VV 
 @VV H_{2} V     @VV k V\\
(\Bbb R,0) @<X+Y << (\Bbb R^{2},0) @>\gamma_{1}>>
(\Bbb R^2,0) @<\gamma_{2}<< (\Bbb R^{2},0) @>\widetilde{f}>>
(\Bbb R, 0)\\
\endCD
\tag 3.1$$
commutes for some diffeomorphism germs $h, H_{1}, K, H_{2}, k$,
where $H_{1}, K, H_{2}$
have the forms in Lemma 3.2,
$f=x+y+$higher term, $\widetilde{f}$ satisfies
$$\widetilde{f}(x,0)=\widetilde{f}(0,x)=x .\tag 3.2 $$
Then from (3.1) we have
$h(x+y)=xA^{2}(x, y^{2})+yB(x, y^{2}),
k\circ f(x, y)=\widetilde{f}(xA\circ \gamma_{2}(x, y), yB^{2}\circ \gamma_{2}(x, y)).$
Thus we have $h(x)=xA^{2}(x, 0)$ and $h(y)=yB(0, y^{2}).$
Therefore $h$ is increasing and odd.
Namely there exists $a\in\Cal E _{1}$ with $a(0)>0$ such that
$$h(t)=ta(t^{2}).$$
Hence $A(u, 0)=\pm\,\root\of{a(u^{2})}$ and $B(0, v)=a(v)$
for any $v\geq 0$.
Moreover from (3.2) we have
$$k\circ f(x, 0)=xA(x^{2}, 0), k\circ f(0, y)=yB^{2}(0, y).\tag 3.3$$
Therefore we have 
$k\circ f(x, 0)=\pm x\,\root\of{a(x^{4})}$ and
$k\circ f(0, y)=ya^{2}(y)$ for $y\geq 0$.

So we now suppose that 
$k\circ f(t, 0)=t\,\root\of{a(t^{4})},
k\circ f(0, t)=ta^{2}(t)$
hold for all $t\in (\Bbb R , 0)$.
The we have 
$$k=t\,\root\of{a(t^4 )}\circ f(t, 0)^{-1}=
ta^2 (t)\circ f(0, t)^{-1},\tag 3.4$$
where $f(t, 0)^{-1}, f(0, t)^{-1}$ respectively denoting the inverse 
function germ of $f(t, 0), f(0, t).$
We can write
$$f(0, t)^{-1}\circ f(t, 0) =tb(t)$$
for some $b\in\Cal E _{1}$ with $b(0)=1.$
Then from (3.4) we have 
$$a(t^{4})=b^{2}(t)a^{4}(tb(t)).
\tag 3.5 $$

Conversely the following holds.
\proclaim{Proposition 3.4} If for any $f\in\Cal E_2$ with the form 
$f(x, y)=x+y+higher\; term$, there exists $a\in\Cal E _1$ with
$a(0)>0$ satisfying the equation (3.5),
then there exists coordinate changes
$h, H_1 , K, H_2 , k$ such that the diagram (3.1) commutes and
$\widetilde{f}=k\circ f\circ H_2 ^{-1}$ satisfying the property (3.2).
\endproclaim
\demo{Proof}
For $f=x+y+higher\; term$, let $a\in\Cal E_1$ with $a(0)>0$ be a solution
of the equation (3.5).
Then from (3.5) we have the equality on the right hand side in (3.4).
So that we define $k(t)$ by (3.4).
Now we put $h(t)=ta(t^2 ).$
Then from the Malgrange preparation theorem
$h(x+y)=x\alpha(x, y^2 )+y\beta(x, y^2 )$ for some
$\alpha, \beta\in\Cal E_{u, v}$ with $\alpha(0)>0.$
By putting $A=\root\of\alpha, B=\beta(u, v)-\beta(0, v)+a(v),$
we have (3.3).
Define $H_1, K, H_2$ by the form in Lemma 3.2.
Then we have the commutativity of the diagram (3.1) and
from (3.3) we see that $\widetilde{f}=k\circ f\circ H_2^{-1}$
satisfies (3.2).
\qed
\enddemo

Therefore the proof of Theorem A is reduced to show,
for any $f(x, y)=x+y+higher\; term,$ the existence of $a\in\Cal E _1$
in Proposition 3.4.
To do this firstly we shall consider
under the formal category and next we shall consider under
$C^{\infty}$ category to kill the flat terms. 

\proclaim{Lemma 3.5} Any germ of type (III, III) with the form in Lemma 3.3
is equivalent to
$$(\Bbb R, 0) @<x+y<< (\Bbb R ^{2}, 0) @>\gamma_{1}>> (\Bbb R ^{2}, 0)
              @<\gamma_{2}<< (\Bbb R^{2}, 0) @>\widetilde{f}>> (\Bbb R, 0),$$
where  $\widetilde{f}(x, 0)=x+flat\; term,\; \widetilde{f}(0, y)=y+flat\; term.$
\endproclaim
\demo{Proof} 
Define the map $j^{\infty} : \Cal E _{1}\to 
\Cal E _{1}/ \Cal M _{1}^{\infty}\cong\Bbb R [[t]]$
by $j^{\infty}(a(t))$=the Talyor expansion of $a(t).$
By direct calculations we see that
for any $\bar{b}(t)\in \Bbb R [[t]]$ with $\bar{b}(0)=1$
there exists 
$\bar{a}(t)\in\Bbb R [[t]]$ satisfying
$$\bar{a}(t^{4})=\bar{b}^{2}(t)\bar{a}^{4}(t\bar{b}(t)),\; \bar{a}(0)>0.
\tag 3.6$$
For any $f =x+y+higher\; term,$ we put $\bar{b}=j^{\infty}(b),$
where $f(0, t)^{-1}\circ f(t, 0)=tb(t).$
Then let $\bar{a}\in\Bbb R [[t]]$ be a solution of (3.6). 
Applying the Borel's theorem there exists $a\in\Cal E _{1}$ such that
$j^{\infty}(a)=\bar{a}$.
Then we put $h(t)=ta(t^{2})$
and moreover
we take a representative element $k$ in
$t\,\root\of{a(t^{4})}\circ f (t, 0)^{-1}+\Cal M_{1}^{\infty}.$
We define a $(\gamma_{1}, \gamma_{2})$-compatible diffeomorphism germ
by the same way as in the proof of Proposition 3.4.
Then by using these coordinate changes 
we have the required result.
\qed
\enddemo 

To kill the flat terms of $\widetilde{f}$ we shall apply Proposition 3.4. 
That is, since $\widetilde{f}(0, t)^{-1}\circ\widetilde{f}(t, 0)=t(1+flat\; function),$
in order to prove Theorem A  
it is sufficient to show 
\proclaim{Proposition 3.6}
For any $b(t)\in\Cal E _1$ of the form
$b(t)=1+flat\; function$
there exists $a\in\Cal E_{1}$ with $a(0)>0$ satisfying (3.5).
\endproclaim

To prove this we need the following Lemma.
Let $\theta : (\Bbb R, 0)\to(\Bbb R, 0)$ be the inverse function germ
of $t\mapsto tb(t),$
where $b(t)=1+flat\; function.$
We see that $\theta(x)=xc(x)$ where $c(x)=1+flat\; function.$
Then from the equation (3.5), we see that
$a(\theta(x)^{4})=(1/c(x))^{2}a(x)^{4}.$
Thus we have   
$$a(x)=c(x)^{\frac{1}{2}}a(\theta(x)^{4})^{\frac{1}{4}}.
\tag 3.7$$

Then by direct calculations inductively we have
\proclaim{Lemma 3.7}
(i)\, In the above equation (3.7), for any non-negative integer $n$, 
$a(x)$ satisfy 
$$a(x)=\left(\prod\limits_{k=0}^{n}c(\sigma^{k}(x))^{1/{4^{k}}}\right)^{1/2}
a(\sigma^{n+1}(x))^{1/{4^{n+1}}},$$
where $\sigma(x)=\theta(x)^{4},$ 
$\sigma^{0}(x)=x, \sigma^{n}=\undersetbrace{n}\to
{\sigma\circ\dots\circ\sigma}\, (n\geq 1).$

(ii)\,$\lim\limits_{n\to\infty}a(\sigma^{n}(x))^{1/4^{n}}=1$
on $(\Bbb R , 0).$
\endproclaim

Hence if $a(x)\in\Cal E _{1}$ satisfying the equation (3.5) exists,
then $a(x)$ must be $a(x)= 
\left(\prod\limits_{n=0}^{\infty}c(\sigma^{n}(x))^{1/{4^{n}}}\right)^{1/2}.$
Actually the following holds:
\proclaim{Proposition 3.8}
$\left(\prod\limits_{n=0}^{\infty}c(\sigma^{n}(x))^{1/{4^{n}}}\right)^{1/2}$
exists and is $C^{\infty}$.
\endproclaim

Therefore proof of Proposition 3.6 is reduced to proof of Proposition 3.8.
To prove Proposition 3.8 we need the following Lemmas.

We recall that
$(\prod\limits_{n=0}^{\infty}f_{n})=\lim\limits_{N\to\infty}f_{0}\cdots f_{N}$
exists and is $C^{\infty}$ if and only if 
$\lim\limits_{N\to\infty}\sum\limits_{k=0}^{N}\log f_{k}$
exists and is $C^{\infty}$,
where $f_{k}\in \Cal E _{1}$ with $f_{k}>0.$

Put $S_{N}(x)=\sum\limits_{n=0}^{N}\log{c}(\sigma ^{n}(x))^{1/{4^{n}}}
=\sum\limits_{n=0}^{N}\frac{1}{4^{n}}\log c(\sigma^{n}(x))$.

\proclaim{Lemma 3.9}
$\left(\prod\limits_{n=0}^{\infty}c(\sigma^{n}(x))^{1/{4^{n}}}\right)^{1/2}$
is uniformly convergent.
\endproclaim
\demo{Proof}
There exists $\varepsilon >0$
such that $|\sigma^{n}(x)|<{|x|}^{3^{n}}<\varepsilon<1,$
which leads to\allowlinebreak
$|\frac{1}{4^{n}}\log c(\sigma^{n}(x))| \allowmathbreak
<\frac{1}{4^{n}}$
for any $n$ on $|x|<\varepsilon$.
Hence $\lim\limits_{N\to \infty}S_{N}(x)$
is uniformly convergent on $|x|<\varepsilon.$ 
\qed
\enddemo

Next we shall show that
$S(x)=\lim\limits_{N\to\infty}S_{N}(x)$ is of class $C^{\infty}.$
Firstly we see that 
$S(x)=\frac{1}{4}S(\sigma(x))+\log c(x).$
Actually from the definition of $S_{N}$ we have
$\frac{1}{4}S_{N}(\sigma(x))=\sum\limits_{n=0}^{N}\frac{1}{4^{n+1}}\log
c(\sigma^{n+1}(x))=S_{N+1}(x)-S_{0}(x).$
Then by $N\to\infty$ we have the required relation.

Therefore if $S(x)$ is of class $C^{p} (p=0, 1, 2, \cdots),$
we inductively have
$$S^{(p)}(x)=\frac{1}{4}S^{(p)}(\sigma(x))\sigma'(x)^{p}
+\lambda_{p}(x), \tag 3.8$$
where $S^{(p)}(x)$=the $p$-th order differential of $S(x)$,
$\lambda_{p}(x)\in\Cal E_{1}$ depend on $p$.
We put $T(x):=S^{(p)}(x),
\alpha(x):=\sigma ' (x)^{p}, \lambda(x):=\lambda_{p}(x).$
Then we inductively have
$$T(x)=\frac{1}{4^n}\alpha(x)\alpha(\sigma(x))\alpha(\sigma^{2}(x))\cdots
\alpha(\sigma^{n-1}(x))T(\sigma^{n}(x))\tag 3.9$$
$$+\sum_{i=0}^{n-1}\frac{1}{4^i}\alpha(x)\alpha(\sigma(x))\cdots
\alpha(\sigma^{i-1}(x))\lambda(\sigma^{i}(x)).
$$

\proclaim{Lemma 3.10}
There exists $\varepsilon >0$ such that
$$\lim_{n\to\infty}\frac{1}{4^n}\alpha(x)\alpha(\sigma(x))\alpha(\sigma^{2}(x))\cdots
\alpha(\sigma^{n-1}(x))T(\sigma^{n}(x))=0$$
on $|x|<\varepsilon$ for any $p$.
\endproclaim
\demo{Proof}
Clearly there exists $\varepsilon \in (0, 1)$ such that
$|\sigma '(x)|<|x|^{2}$ on $|x|<\varepsilon.$
Thus we see that
$$|\frac{1}{4^{n}}\sigma'(x)^{p}\sigma'(\sigma(x))^{p}\cdots
\sigma'(\sigma^{n-1}(x))^{p}T(\sigma^{n}(x))|
<\frac{1}{4^{n}}|x|^{2pn}|T(\sigma^{n}(x))|$$
on $|x|<\varepsilon$ for any $p$.
Then the right hand side is convergent to $0$ on $|x|<\varepsilon$
if $n\to\infty$.
\qed
\enddemo

\proclaim{Lemma 3.11}
There exists $\varepsilon >0$ such that for any $\alpha (x)\in\Cal E_{1}$
with $\alpha(0)=\alpha '(0)=0$ and 
for any $\lambda(x)$ which is $C^{1}$ class function germ on $(\Bbb R, 0)$
with $\lambda(0)=\lambda'(0)=0$ the followings hold:
Let $u_{n}(x)=\alpha(x)\alpha(\sigma(x))\alpha(\sigma^{2}(x))\cdots
\alpha(\sigma^{n}(x))\lambda(\sigma^{n+1}(x)).$

{\text{(i)}}$\lim\limits_{N\to\infty}\sum\limits_{n=0}^{N}u_{n}(x)$ 
is pointwise convergent on $|x|<\varepsilon.$

{\text{(ii)}}$\lim\limits_{N\to\infty}\sum\limits_{n=0}^{N}u_{n}(x)$
is of class $C^1.$ 
\endproclaim
\demo{Proof of (i)}
We can take an $\varepsilon_0\in (0, \frac{1}{2})$ such that
for any $\alpha, \lambda$ in this Lemma there exists a positive number
$M_{\alpha, \lambda}\in\Bbb R$ such that the followings hold for any
$|x|<\varepsilon_0$:
$$|\alpha(x)|<M_{\alpha, \lambda},\,
|\alpha '(x)|<M_{\alpha, \lambda},\,
|\lambda(x)|\leq M_{\alpha, \lambda}|x|,\,
|\lambda'(x)|\leq M_{\alpha, \lambda}|x|.$$
Take $\varepsilon_{\sigma}\in (0, \frac{1}{2})$ such that
$|\sigma|<|x|^3$ (hence $|\sigma ^n |<|x|^{3^n}$ for any $n$)
on $|x|<\varepsilon_{\sigma}.$
Let $\varepsilon_1 = \min \{\varepsilon_0 , \varepsilon_\sigma \}.$
Then we see $|u_n | < M_{\alpha, \lambda}^{n+2}|\sigma^{n+1}(x)|
< M_{\alpha, \lambda}^{n+2}|x|^{3^{n+1}}$ on $|x|<\varepsilon_1 . $
We easily verify that the series $\sum M_{\alpha, \lambda}^{n+2}|x|^{3^{n+1}}$ converges pointwise on $|x|<\varepsilon _1$
and hence $\sum u_n$ is also.\qed
\enddemo
\demo{Proof of (ii)}
Since $u_n (x)$ is of class $C^1$, it is sufficient to show that there exists $\varepsilon > 0$ independent of $\lambda$ such that
the series $\sum u_n ' (x)$ converges uniformly on $|x|<\varepsilon.$
Take $\varepsilon _2 >0$ such that $|(\sigma ^i )' |<1$
for any $i$ on $|x|<\varepsilon _2 .$
Let $\varepsilon = \min\{\varepsilon _1 , \varepsilon _2\},$ where 
$\varepsilon _1$ is defined just as before.
Then we have $|u_n ' (x)|<(n+2)M_{\alpha, \lambda}^{n+2}|x|^{3^{n+1}}$
on $|x|<\varepsilon.$ Hence we easily verify that
$\sum u_n ' (x)$ converges pointwise.

Moreover we have $|u_n ' (x)|<(n+2)M_{\alpha, \lambda}^{n+2}
(1/2^{3^{n+1}}) =: K_n$
on $|x|<\varepsilon.$ Since $\root{n}\of{K_n} \to 0,$
$\sum K_n$ converges.
Therefore $\sum u_n ' (x)$ converges uniformly on $|x|<\varepsilon.$
\qed
\enddemo

Now we are ready to show that $S(x)$ is of class $C^{\infty}.$

\demo{Proof of Proposition 3.8}
We recall the relation (3.8) of $S(x)$.
In (3.8), if $S(x)$ is of class $C^p$, by direct calculations we see
$\lambda_0 = \log c, \lambda_1 = \frac{c'}{c}, \lambda_p =
\frac{1}{4}S^{(p-1)}(\sigma)\cdot(p-1)(\sigma ')^{p-2}\sigma'' 
+\lambda_{p-1}'$
=sums and products of $\lambda_1^{(p-1)}, S'(\sigma), S''(\sigma),\dots,
S^{(p-1)}(\sigma), \sigma, \sigma', \sigma'',\dots, \sigma^{(p)}$
for $p\geq 2.$
Since $c, \sigma$ are of class $C^\infty$, if $S(x)$ is of class $C^p$,
then $\lambda_p$ is of class $C^{1}.$
Then by putting $\lambda=\lambda_p , \alpha=(\sigma ')^p$, from
the relation (3.9), Lemma 3.10 and Lemma 3.11 we can inductively
prove that $S^{(p)}=\sum\limits_{n=0}^{\infty}\frac{1}{4^{n+1}}u_n
+\lambda_p$ is of class $C^1$ for any $p.$\qed
\enddemo

\head{4. Topological rigidity for type (III, III)}
\endhead

In this section we shall prove the following which is a part
of Theorem B:

\proclaim{Proposition 4.1}
If two diagrams of type (III, III) are topologically equivalent,
then they are $C^{\infty}$-equivalent.
\endproclaim

One of the essential tool to prove theorem B is the following Theorem
of rigidity of webs.
Let $\Cal W =(\Cal F _1, \dots ,\Cal F _d )$ be a configuration of d
foliations $\Cal F _j , j=1, \dots , d (d\geq 2)$ in a domain $U$ 
of $\Bbb R ^2$ with codimension 1.
Define the set $S(\Cal W )$ by $S(\Cal W)$ :=
the set of points at which $\Cal F _i$ and $\Cal F _j$ ($1\leq i <j\leq d$)
are non-transversal.
We call $\Cal W$ a {\it non-singular d-web\/} or briefly {\it d-web}
(resp.  {\it singular d-web\/})
on $U$ of codim 1 if $S(\Cal W) =\phi$ (resp. $S(\Cal W )\not=\phi$) and 
call $S(\Cal W)$ the {\it singular set of $\Cal W .$} 
Therefore locally a (singular) d-web is a configuration 
$\Cal W =(f_1,\dots , f_d)$  of d 
submersion germs $f_{j} : (\Bbb R^{2}, 0)\to(\Bbb R, 0).$
Then the singular set $S(\Cal W)$ is given by
$S(\Cal W) =\{p\in(\Bbb R ^2 , 0)\, |\, df_i \wedge df_j (p)=0,\, 1\leq i < j \leq d\}.$ 
 
\proclaim{Theorem 4.2}(Rigidity of webs[4])
Let $\Cal W =(f_{1}, f_{2}, f_{3}), 
\Cal W' =(f'_{1}, f'_{2}, f'_{3})$ 
be two 3-webs on $(\Bbb R^{2}, 0).$
If there exist homeomorphism germs
$\Phi : (\Bbb R ^{2}, 0)\to(\Bbb R ^{2}, 0),
 h_{j}:(\Bbb R , 0)\to (\Bbb R, 0)$
such that $h_{j}\circ f_{j}=f'_{j}\circ\Phi$
for any $j=1,2,3$,
then $\Phi, h_j (j=1, 2, 3)$ are $C^{\infty}$-diffeomorphism germs.
\endproclaim

As a corollary of
Theorem 4.2 we see that rigidity of webs holds for every d-webs with d$\geq 3$.    

Suppose that two diagrams of type (III, III) are topologically equivalent.
That is, from Theorem A, we suppose that the following diagram commute:
$$\CD
(\Bbb R,0) @<x+y << (\Bbb R^{2},0) @>\gamma_{1}>>
(\Bbb R^2,0) @<\gamma_{2}<< (\Bbb R^{2},0) @>f>>
(\Bbb R, 0)\\
 @V h_{1} VV         @V H_{1} VV        @V K VV 
 @VV H_{2} V     @VV h_{2} V\\
(\Bbb R,0) @<X+Y << (\Bbb R^{2},0) @>\gamma_{1}>>
(\Bbb R^2,0) @<\gamma_{2}<< (\Bbb R^{2},0) @>\widetilde{f}>>
(\Bbb R, 0)\\
\endCD
\tag 4.1$$
where $f,\, \widetilde{f}\in\Cal E_{2}$ with
$f(x, 0)=f(0, x)=x,\,   \widetilde{f}(\widetilde{x}, 0)
=\widetilde{f}(0, \widetilde{x})=\widetilde{x}$
and
$h_{i},\,  H_{i} (i=1, 2),\, K$ are homeomorphism germs.
We put $H_{1}=(X, Y),\, H_{2}=(\widetilde{X}, \widetilde{Y}),\,
K=(U, V).$

We remark that there exists a 4-web in the domain
$\{(u, v)\in (\Bbb R ^{2}, 0) | u>0\; {\text{and}}\; v>0\}$
which are constructed by $\gamma_{i}({f_{i}}^{-1}(t)), i=1, 2.$
Therefore from Theorem 4.2 we see that
$K$ restricted to $\{(u, v)\in (\Bbb R ^{2}, 0) | u>0 \;{\text{and}}\; v>0\}$
is  a $C^{\infty}$-diffeomorphism germ.

\proclaim{Lemma 4.3}
We have the followings:

(i) $h_{1}, h_{2}$ are $C^{\infty}$-diffeomorphism germs and 
moreover  $h_{1}=h_{2}=id_{\Bbb R}$ 

(ii) $H_{1}=id_{\Bbb R^{2}}$
  
(iii) $H_{2}$ is a $C^{\infty}$-diffeomorphism germ and moreover
$H_{2}$ restricted to $\{(x,y)\in(\Bbb R^{2},0) | 
\mathbreak 
y\geq0\;{\text{or}}\; x=0\}$
is the identity map germ.
\endproclaim
\demo{Proof of (i)}
We take a point $p\in\{x+y=0\}\cap\{x>0, y<0\}$ such that
$p$ is sufficiently near 0.
Then we consider a non-singular $C^{\infty}$-curve $T$ at $p$
such that $T$ and each lines $x+y=t (t\in(\Bbb R , 0))$ are
transversal. Then $x+y$ restricted to $T$ is a $C^{\infty}$-diffeomorphism
germ.
Since $\gamma_1|_{\{x>0, y<0\}}$ and $K|_{\{u>0, v>0\}}$ are $C^{\infty}$-
diffeomorphism germs, from (4.1) we see that
$H_1|_{\{x>0, y<0\}}$ is a $C^{\infty}$-diffeomorphism germ.
From (4.1) we have $h_1\circ (x+y)|_T  =(X+Y)\circ H_1 |_T $,
namely $h_1$ is a $C^{\infty}$-diffeomorphism germ.
For $h_2$, similarly take a point $p\in f^{-1}(0)\cap\{x<0, y>0\}$ 
which is sufficiently near 0 and then consider
a curve $T$ at $p$ such that $T$ and each fibres of $f$ are transversal.
Then by the similar argument we see that $h_2$ is a $C^{\infty}$-
diffeomorphism germ.

Hence there exist $A, B\in\Cal E_{u, v}$ such that
$h_{1}(x+y)=A(x, y^{2})+yB(x, y^{2}).$
Similarly there exist $\widetilde{A}, \widetilde{B}\in\Cal E_{u, v}$
such that
$h_{2}\circ f(x, y)= \widetilde{A}(x^{2}, y)+x\widetilde{B}(x^{2}, y).$
We remark that in (4.1)
$H_{1}, K, H_{2}$ preserve both the horizontal-axis and the vertical-axis
in each spaces.
From (4.1) we have $U(x, y^2 )=A(x, y^2 ), V(x, y^2 )=y^2 B(x, y^2 )^2$ 
and similarly
$U(x^2, 0)=x^2 \widetilde{B}(x^2 , 0)^2 , V(0, y)=\widetilde{A}(0, y).$
Thus we have
$h_{1}(y)=yB(0, y^{2})$ and $h_{1}(x)=x\widetilde{B}(x, 0)^{2}$
for $ x\geq 0$.
Similarly
$h_{2}(x)=x\widetilde{B}(x^{2}, 0)$ and $h_{2}(y)=yB(0, y)^{2}$
for $ y\geq 0$.
This means that $h_{1}, h_{2}$ are increasing and odd.
Therefore there exist $a_{i}\in\Cal E_{1}$ with $a_{i}(0)>0$
such that $h_{i}(t)=ta_{i}(t^{2})$, $i=1, 2.$
Hence $a_{i}(t)^{4}=a_{i}(t^{4})$ for $t\geq0$,
$i=1,2.$
Thus we see $a_{i}(t)\equiv1$ for $t\geq0$,
namely $h_{i}(t)=t$ on $(\Bbb R, 0).$
\qed 
\enddemo
\demo{Proof of (ii)}From (i) it is clear.
\qed
\enddemo
\demo{Proof of (iii)}
For $f=x+y+$higher term $\in \Cal E_{2}$ consider the map germ
$J_{f} : (\Bbb R^{2}, 0)\to(\Bbb R^{2},0)$ defined by
$J_{f}(x,y)=(f(x, y), f(-x, y)).$
Clearly $J_{f}$ is a $C^{\infty}$-diffeomorphism germ.
Since $H_{2}(-x, y)=(-\widetilde{X}(x, y), \widetilde{Y}(x, y))$
and $h_{2}=id$, we have 
$J_{f}=J_{\widetilde{f}}\circ H_{2}.$
Namely $H_{2}$ is a $C^{\infty}$-diffeomorphism germ.
The latter part of (iii) is clear by (i) and (ii).
\qed
\enddemo

From Lemma 4.3 (iii) and the $\gamma_2$-compatibility,
we can put $H_2=(xA(x^2, y), yB(x^2, y))$ for some
$A, B\in\Cal E _{u, v}$ with $A=1$ on $\{u\geq 0\;{\text{and}}\; v\geq 0\},
B=1$ on $\{u\geq 0\;{\text{and}}\; v\geq 0\}\cup\{u=0\}.$
Then we see that $K=(uA^2, vB)$ on $\{u\geq0\}.$
On the other hand $K(u, v)=(u, v)$ on $\{v\geq0\}.$
Now we need the following Proposition.
\proclaim{Proposition 4.4}
Let $h\in\Cal E _{u, v}$ with $h=1$ on $\{u\geq 0\; {\text{and}}\;
v\geq 0\}.$ Then there exists
$\widetilde{h}\in\Cal E _{u, v}$
such that
$\widetilde{h}=h$ on $\{u\geq0\}$ and
$\widetilde{h}=1$ on $\{u\leq0\;{\text{and}}\; v\geq0\}.$
\endproclaim

If this Proposition 4.4 is proved, the proof of Proposition 4.1
is finished by the following.
\demo{Proof of Proposition 4.1}
For the above $A, B\in\Cal E_{u, v}$ we take its extension
$\widetilde{A}, \widetilde{B}$ with the same property in Proposition 4.4.
Then define the $C^{\infty}$-diffeomorphism germ 
$\widetilde{K} : (\Bbb R^2, 0)\to (\Bbb R^2, 0)$ by
$\widetilde{K}=(u\widetilde{A}(u, v)^2, v\widetilde{B}(u, v)).$
Then we see that $\widetilde{K}\circ\gamma_1 =\gamma_1$
and $\widetilde{K}\circ\gamma_2 = \gamma_2 \circ H_2 .$
This completes the proof.\qed
\enddemo

In order to prove Proposition 4.4 
we apply the theory of 
{\it{regularly situated sets}} by 
$\L$ojasiewicz as follows (See [11], [12]p.12):
Let $\Omega$ be an open set in $\Bbb R ^{n}$ and let
$X$ be a closed subset of $\Omega$. Let $\Cal E (X)$ be the ring of
all Whitney functions of infinite order on $X.$ 

Let $X, Y$ be closed subsets of $\Omega.$
Let $\delta : \Cal E(X\cup Y)\to\Cal E (X)\oplus\Cal E (Y)$ 
be the diagonal mapping defined by
$$\delta(F)=(F|_{X}, F|_{Y}).$$
Let $\pi : \Cal E (X)\oplus\Cal E (Y)\to\Cal E (X\cap Y)$ be the mapping
defined by
$$\pi(F, G)=F|_{X}-G|_{Y}.$$

\pagebreak
\definition{Definition 4.5}([11, 12])
Two closed subsets $X, Y$ of an open set $\Omega$ are said to be
{\it{regularly situated}} if the sequence
$$0 @>>> \Cal E(X\cup Y) @>\delta>> \Cal E (X)\oplus\Cal E (Y)
@>\pi>> \Cal E(X\cap Y) @>>>0
$$
is exact.
\enddefinition

\proclaim{Theorem 4.6}($\L$ojasiewicz)
Given $X, Y$ closed in an open set $\Omega$ a necessary and sufficient
condition that they are regularly situated is the following:
Either $X\cap Y=\phi$ or

$(\Lambda)$ Given any pair of compact sets $K\subset X, L\subset Y,$
there exists a pair of constants 
$C>0$ and $\alpha>0$ such that, for every $x\in K$, one has
$$d(x, L)\geq C\, d(x, X\cap Y)^{\alpha},$$
where $d$ denoting the euclidean distance in $\Bbb R^n.$
\endproclaim

We apply Theorem of $\L$ojasiewicz for the following particular case.
Consider closed subsets $X, Y$ of $\Bbb R^{2}$ defined by 
$$X=\{r{e^{i\theta}} | 0\leq r\leq\varepsilon, 
-\frac{\pi}{4}\leq\theta\leq\frac{2}{3}\pi\},Y= \{r{e^{i\theta}} | 0\leq r\leq\varepsilon, \frac{5}{6}\pi\leq\theta\leq\frac{7}{4}\pi\},$$
where $\varepsilon$ is a fixed positive number.

\proclaim{Lemma 4.7}$X, Y$ are regularly situated.
\endproclaim
\demo{Proof}
We can easily see $d(p, Y)=d(p, X\cap Y)$ for every 
$p\in \{r{e^{i\theta}} | 0\leq r\leq\varepsilon, 
-\frac{\pi}{4}\leq\theta\leq\frac{\pi}{3}\}$
and also see that
$d(p, Y)=d(p, 0)\sin (\frac{5}{6}\pi - \arg p),\, 
d(p, X\cap Y)=d(p, 0)$ for every 
$p\in \{r{e^{i\theta}} | 0\leq r\leq\varepsilon, 
\frac{\pi}{3}\leq\theta\leq\frac{2}{3}\pi\}.$
Thus we have
$d(p, Y)\geq \frac{1}{2}d(p, X\cap Y)$ for any $p\in X.$
Then for $X, Y$ obviously the condition $(\Lambda)$ holds.
\qed\enddemo

Clearly to prove Proposition 4.4 it is enough to show the following.
\proclaim{Lemma 4.8}
Let $h : (\Bbb R^{2}, 0)\to(\Bbb R, 0)$ be a $C^{\infty}$-function germ
such that
$h=0$ on $\{(x, y) | x\geq0\; {\text{and}}\; y\leq0\}$.
Then there exists a $C^{\infty}$-function germ
$\widetilde{h} : (\Bbb R^{2}, 0)\to(\Bbb R, 0)$ such that
$\widetilde{h}=h$ on $\{x\geq0\}$ and
$\widetilde{h}=0$ on $\{y\leq0\}.$
\endproclaim
\demo{Proof}
We take a sufficiently small $\varepsilon >0.$
Then we can define $f\in\Cal E(X)$ by $h|_{X}$ and 
define $g\in\Cal E(Y)$ by $0|_{Y}.$
Obviously $(f, g)\in$ Ker $\pi.$
Then by Lemma 4.7 there exists $F\in\Cal E(X\cup Y)$
such that
$\delta (F)=(f, g).$
By the definition of Whitney function there exists a $C^{\infty}$-function
$\widetilde{F}$ defined in a neighbourhood of $X\cup Y$ in $\Bbb R^2$
such that $j^{\infty}\widetilde{F}(p)=F(p)$ in $X\cup Y.$
Put $\widetilde{h}:=$the germ of $\widetilde{F}$ at 0,
which is the required  $C^{\infty}$-function germ.
\qed
\enddemo

\head{5. Rigidity for each types except (III, III)}
\endhead
In this section we shall complete proof of Theorem B.
Firstly we recall the following normal forms which have been given 
in [3, 10, 14]. 
\proclaim{Theorem 5.1}The normal forms of each generic
types except $(III, III)$ for the diagrams 
$(\Bbb R, 0) @<f_1<< (\Bbb R ^{2}, 0) @>\gamma_{1}>> (\Bbb R ^{2}, 0)
              @<\gamma_{2}<< (\Bbb R^{2}, 0) @>f_2>> (\Bbb R, 0)$
are the followings:

$(I, I)^{0}\;\; f_{1}=y_{1},\;\; \gamma_{1}=(x_1, y_1);$

$\;\;\;\;\;\;\;\;\;\;\;\; f_{2}=x_{2},\;\; \gamma_{2}=(x_2, y_2).$

\smallpagebreak

$(I, I)^{1}\;\; f_{1}=y_{1},\;\; \gamma_{1}=(x_1, y_1);$

$\;\;\;\;\;\;\;\;\;\;\;\;  f_{2}=x_{2}^{2}+y_{2},\;\; \gamma_{2}=(x_2, y_2).$

\smallpagebreak

$(I, I)^{2}\;\; f_{1}=y_{1},\;\; \gamma_{1}=(x_{1}, y_{1});$

$\;\;\;\;\;\;\;\;\;\;\;\;  f_{2}=x_{2}^{3}+x_{2}y_{2}+y_{2},\;\;
\gamma_{2}=(x_2, y_2).$

\smallpagebreak

$(II, I)\;\; f_{1}=x_{1}^2\pm y_{1}^2,\;\; \gamma_{1}=(x_1, y_1);$

$\;\;\;\;\;\;\;\;\;\;\;\;  f_{2}=x_{2},\;\; \gamma_{2}=(x_{2}, y_{2}).$

\smallpagebreak

$(III, I)^{0}\;\; f_{1}=x_{1}+y_{1},\;\; \gamma_{1}=(x_1, y_{1}^{2});$

$\,\;\;\;\;\;\;\;\;\;\;\;\;\;\;\; f_{2}=x_{2}+\theta(x_{2}, y_{2}),\;\; 
                   \gamma_{2}=(x_2, y_2),$

$\,\;\;\;\;\;\;\;\;\;\;\;\;\;\;\; \text{where}\; \theta\in\Cal M _{x_{2}, y_{2}} 
\;\text{with}\; \theta(x_{2}, 0)=0.$

\smallpagebreak

$(III, I)^{1}\;\; f_{1}=x_{1}+y_{1},\;\; \gamma_{1}=(x_1, y_{1}^{2});$

$\,\;\;\;\;\;\;\;\;\;\;\;\;\;\;\; f_{2}=x_{2}^{2}+\theta(x_{2}, y_{2}),\;\;
 \gamma_{2}=(x_{2}, y_{2}),$

$\,\;\;\;\;\;\;\;\;\;\;\;\;\;\;\; \text{where}\; \theta\in\Cal M _{x_{2}, y_{2}}
\;\text{with}\;  \theta(x_{2}, 0)=0,
\frac{\partial \theta}{\partial y_{2}}(0)\not= 0.$

\smallpagebreak

$(IV, I)\;\; f_{1}= x_{1}^{2}+y_{1},\;\; \gamma_{1}=(x_1, y_{1}^{2});$

$\;\;\;\;\;\;\;\;\;\;\;\;  f_{2}=x_{2}+\theta(x_{2}, y_{2}),\;\; \gamma_{2}=(x_2, y_2),$

$\;\;\;\;\;\;\;\;\;\;\;\; \text{where}\; \theta\in\Cal M _{x_{2}, y_{2}} 
\;\text{with}\; \theta(x_{2}, 0)=0, 
\frac{\partial \theta}{\partial y_{2}}(0)\not= 0.$

\smallpagebreak

$(V, I)\;\; f_{1}=x_{1}+x_{1}y_{1}+y_{1}^{3},\;\; 
\gamma_{1}=(x_1, y_{1}^{2});$

$\;\;\;\;\;\;\;\;\;\;\; f_{2}=x_{2}+\theta(x_{2}, y_{2}),\;\;
 \gamma_{2}=(x_2, y_2),$

$\;\;\;\;\;\;\;\;\;\;\; \text{where}\; \theta\in\Cal M _{x_{2}, y_{2}}
\;\text{with}\; \theta(x_{2}, 0)=0,
\frac{\partial \theta}{\partial y_{2}}(0)\not= 0,
\frac{\partial \theta}{\partial y_{2}}(0)\not= 3.$

\smallpagebreak

$(VI, I)\;\; f_{1}=y_{1}+\alpha\circ \gamma_{1},\;\; 
\gamma_{1}=(x_1, y_{1}^{3}+x_{1}y_{1});$

$\;\;\;\;\;\;\;\;\;\;\;\;\; f_{2}=x_{2}+\theta(x_{2}, y_{2}),\;\; \gamma_{2}=(x_2, y_2),$

$\;\;\;\;\;\;\;\;\;\;\;\;\; \text{where}\; \alpha\in \Cal M _{u, v},\,
\theta\in\Cal M _{x_{2}, y_{2}} \;\text{with}\; \theta(x_{2}, 0)=0.$

\endproclaim

For the normal forms with functional moduli except type $(III, III)$, that is,
for the normal forms of $(III, I)^{0}, (III, I)^{1}, (IV, I),
(V, I)$ and $(VI, I)$ 
we may skip $\gamma_2$ to consider the topological rigidity problem
because $\gamma_2 =id.$

\proclaim{5.1 Proof of rigidity for (III, I)$^{0}$, (IV, I), (V, I)}
\endproclaim
We uniformly prove for these types.
In this proof we always suppose that $\gamma = (x, y^2 ) : 
(\Bbb R^2, 0)\to (\Bbb R^2, 0).$
Assume that two normal forms of the same type(i.e. (III, I)$^{0}$,
(IV, I), (V, I)) are topologically equivalent:  
$$\CD
(\Bbb R,0) @<g << (\Bbb R^{2},0) @>\gamma>>
(\Bbb R^2,0) @>f>>  (\Bbb R, 0)\\
 @V h VV         @V H VV        @V K VV      @VV k V\\
(\Bbb R,0) @<g << (\Bbb R^{2},0) @>\gamma>>
(\Bbb R^2,0) @>\widetilde{f}>> (\Bbb R, 0)\\
\endCD
\tag 5.1$$
where $h, H, K, k$ are homeomorphism germs.

\proclaim{Step 1} $h, k$ are $C^{\infty}$-diffeomorphism germs.
\endproclaim

\demo{Case $(III, I)^{0}$}
The diagrams of this type have a 3-web structure on $\{v>0\}.$
By the rigidity of webs, $K|_{\{v>0\}}$ is a $C^{\infty}$-diffeomorphism
germ. 
Since $\gamma|_{\{y>0\}}$ is a $C^{\infty}$-diffeomorphism germ,
from (5.1) we see that $H|_{\{y>0\}}$ is a $C^{\infty}$-diffeomorphism
germ. Take a point $p\in g^{-1}(0)\cap \{y>0\}$ such that $p$ is sufficiently
near 0. Then we can consider a non-singular curve $T$ at $p$ 
such that $T$ and each fibres of $g$ are transversal.
Namely $g|_T $ gives a $C^{\infty}$-diffeomorphism germ.
Hence $g$ restricted to $H(T)$ is a $C^{\infty}$-diffeomorphism germ.
Then from (5.1) $h\circ g|_{T}=g|_{H(T)}\circ H|_{T},$
namely $h$ is a $C^{\infty}$-diffeomorphism germ.

Similarly we take a point $q\in f^{-1}(0)\cap\{v>0\}$ such that
$q$ is sufficiently near 0 and
then consider a non-singular curve $S$ at $q$ such that
$S$ and each fibres of $f$ are transversal.
Then we see that $k=\widetilde{f}|_{K(S)}\circ K|_{S}\circ f|_{S}^{-1}$
is a $C^{\infty}$-diffeomorphism germ.
\enddemo
\demo{Case $(IV, I)$} In this case there is a singular 3-web on the 
domain $\{v>0\}$ and its singular set is the $v$-axis$(v>0).$
That is, there is a 3-web on the domain 
$\{(u, v) | u\not=0, v>0\}.$
By the rigidity of webs, $K$ restricted to $\{u\not=0, v>0\}$ is 
a $C^{\infty}$-diffeomorphism germ.
Hence we see that $H$ restricted to $\{x>0 ,  \, y<0\}$
is a $C^{\infty}$-diffeomorphism germ.
Take a point $p\in g^{-1}(0)\cap\{x>0, y<0\}$ such that $p$ is sufficiently
near 0. Then taking a curve $T$ at $p$ with the same property as the case $(III, I)^{0}$
we see that $h$ is a $C^{\infty}$-diffeomorphism germ.

On the other hand by the generic condition, since $f^{-1}(0)$ 
and $v$-axis$(v>0)$ are transversal 
we can take a point $q\in f^{-1}(0)\cap\{u\not=0, v>0\}.$
Then we can take a curve $S$ at $q$ in $\{u<0, v>0\}$ or $\{u>0, v>0\}$
which has the same property as in the case $(III, I)^{0}.$
Hence $k$ is a $C^{\infty}$-diffeomorphism germ.
\enddemo
\demo{Case $(V, I)$}
In this case there is a singular 3-web  on $\{v>0\}$
such that 3 foliations are given by the level curves $\gamma(g^{-1}(t))$ and 
$f^{-1}(t).$ That is, the singular 3-web $\Cal W$ is given by
$$f_1 =u+(u+v)\,\root\of{v},\, f_2=u-(u+v)\,\root\of{v},\,
f_3 =f=u+bv+\theta (u, v)$$
for $v>0$, where $b\not=0, \theta\in\Cal M_{u, v}^{2}.$
We decompose the singular set $S(\Cal W)$ by
$S(\Cal W)=S_{1, 2}\cup S_{2,3}\cup S_{3,1},$
where $S_{i, j}=\{(u, v) | v>0, det\, J_{(f_i , f_j)}(u, v)=0\}
(i\not=j)$,
$J_{(f_i , f_j)}(u, v)$ denoting the Jacobi matrix of $(f_i , f_j )$
at $(u, v).$
We have $S_{1, 2}=\{u+3v=0, v>0\}.$
We see that the tangent directions of
the level curves of $f_1$ and $f_2$ are vertical
at every points in $S_{1, 2}.$

To show the Step 1, we shall apply the same technique as the case (III, I)$^0 .$ 
So that we shall avoid the set $S_{2, 3}\cup S_{3, 1}$ as follows.
Let $m\in\Bbb R .$
Then we have 
$det\, J_{(f_3 , f_1)}(mv, v)=(1+\frac{\partial\theta}{\partial u}(mv, v))
\frac{m+3}{2}\,\root\of{v}
-(b+\frac{\partial\theta}{\partial v}(mv, v))(1+\,\root\of{v}).$
Since $\lim\limits_{v\to 0} det\, J_{(f_3 , f_1)}(mv, v)=-b\not=0,$
there exists $\delta_1 (m)>0$ such that
$det\, J_{(f_3 , f_1)}(mv, v)\not=0$ for $\delta_1 (m)>v>0.$
We can choose $\delta_1 (m)$ such that $\delta_1 (m)$ is a continuous
function with respect to $m.$
Similarly there exists $\delta_2 (m)>0$ such that
$det\, J_{(f_3 , f_2)}(mv, v)\not=0$ for $\delta_2 (m)>v>0.$
Now for $f$ we take a sufficiently large $l\in\Bbb R$ such that
$f^{-1}(0)$ and $S_{1, 2}$ lie in the set
$\{(mv, v) | -l<m<l, v>0\}.$
Put $\delta=\min \{\delta_1 (m), \delta_2 (m)\}$ for
$-l\leq m \leq l$ and consider the domain
$D_f =\{(mv, v) | -l<m<l, 0<v<\delta\}.$
Then we see that there is a non-singular 3-web on $D_f -S_{1, 2}.$
By rigidity of webs , we see that $K$ restricted to  $D_f -S_{1, 2}$
is a $C^{\infty}$-diffeomorphism germ.
Therefore $H$ restricted to $\gamma^{-1}(D_f -S_{1, 2})\cap\{y>0\}$
is a $C^{\infty}$-diffeomorphism germ.
Take a point $p\in g^{-1}(0)\cap\{y>0\}$ such that
$p$ is sufficiently near 0.
We can take a curve $T$ at $p$ in 
$\gamma^{-1}(D_f -S_{1, 2})\cap\{y>0\}$ such that
$T$ and each fibres of $g$ around 0 are transversal.
Then by the same argument as the case (III, I)$^0$,
we see that $h$ is a $C^{\infty}$-diffeomorphism germ.

Since $f^{-1}(0)$ and $S_{1, 2}$ are transversal by the generic condition,
we can take a point $q\in f^{-1}(0)\cap (D_f -S_{1, 2})$ and moreover 
we can take a curve $S$ at $q$ in $D_f -S_{1, 2}$
such that $S$ and each fibres of $f$ around 0 are transversal.
Then by the same argument as the case (III, I)$^0$,
we see that $k$ is a $C^{\infty}$-diffeomorphism germ.
\enddemo
\proclaim{Step 2}
$H$ is a $C^{\infty}$-diffeomorphism germ.
\endproclaim
For the normal forms $(g, \gamma, f)$ of each cases
we define a $C^{\infty}$-map germ 
$\Phi _f : (\Bbb R^2, 0)\to (\Bbb R^2, 0)$ by
$\Phi _f =(g, f\circ \gamma).$
Then we have $h\times k\circ\Phi_f = \Phi_{\widetilde{f}}\circ H.$
\demo{Cases $(III, I)^{0}, (IV, I)$} we see that $\Phi_f$ is a 
$C^{\infty}$-diffeomorphism germ. Hence from Step 1, $H$ is a $C^{\infty}$
-diffeomorphism germ.
\enddemo
\demo{Case $(V, I)$} From a generic condition
``the tangent cone of $\gamma(g^{-1}(0))$ and $f^{-1}(0)$ are transversal" 
(i.e.$\frac{\partial f}{\partial v}(0)\not=0$),
by direct caluculations we see that $\Phi_f$ is a fold map germ.
Thus by Step 1 we see that $H$ is a $C^{\infty}$-diffeomorphism germ
because if $\gamma\circ H_1 = H_2\circ\gamma$ , where $H_1$ is a homeomorphism germ and $H_2 =(U, V)$ is a $C^{\infty}$-diffeomorphism germ, 
then we have $H_1=(U(x, y^{2}), y\,\root\of {V_{1}(x, y^{2})})$
for some $V_{1}\in\Cal E _2$ with $V_1 (0)\not=0.$
\enddemo
\proclaim{Step 3}
A $C^{\infty}$-extension of $K$ restricted to the web domain.
\endproclaim
We shall construct a $C^{\infty}$-diffeomorphism germ $\widetilde K$
preserving the commutativity (5.1)
even if we replace $K$ with $\widetilde{K}$ as follows:
By the coordinate change $(f, v)$ in the source of $f$,
we may suppose $f=\pi_{1}$, where $\pi_{1}$ is the first projection
$\pi_{1}(u, v)=u.$ Similarly $\widetilde{f}=\pi_{1}.$
Hence it is enough to construct a $C^{\infty}$-diffeomorphism germ
$\widehat{K} : (\Bbb R^2 , 0)\to (\Bbb R^2, 0)$ such that
the following diagram commutes
$$\CD
(\Bbb R^{2},0) @>\delta>>
(\Bbb R^2,0) @>\pi_1 >>  (\Bbb R, 0)\\
 @V H VV        @V \widehat{K} VV      @VV k V\\
(\Bbb R^{2},0) @>\widetilde{\delta}>>
(\Bbb R^2,0) @>\pi_1>> (\Bbb R, 0)\\
\endCD
$$
where $\delta=(f, v)\circ\gamma, \widetilde{\delta}=(\widetilde{f}, v)
\circ\gamma.$
Let $H=(X, Y).$ 
Then from Step 2 and by the Malgrange preparation
theorem there exists $A, B\in \Cal E_{u, v}$ such that
$Y=A\circ\delta+yB\circ\delta.$
From (5.1)
we see that $A\circ\delta=0$ in $(\Bbb R ^2, 0).$
Then put $\widehat{K}(u, v)=(k(u), vB(u, v)^{2}).$
From (5.1) and Steps 1, 2
we see that $\widehat{K}$ is the required $C^{\infty}$-diffeomorphism germ.

This completes the proof of rigidities for $(III, I)^{0}, (IV, I), (V, I).$
\qed

\proclaim{5.2 Proof of rigidity for (III, I)$^{1}$} 
\endproclaim
The following was shown in [4].
\proclaim{Lemma 5.2}
Let $(\Bbb R, 0) @<f<< (\Bbb R ^{2}, 0) @>\gamma>> (\Bbb R ^{2}, 0)$
be of type (III). Then only using $C^{\infty}$-coordinate changes
in the source and the target of $\gamma$, without any coordinate
changes in the target of $f$, the normal form of (III) is given by
$f=x+y, \gamma=(x, y^2).$
\endproclaim

We suppose that $\gamma =(x, y^2) : (\Bbb R^2, 0)\to (\Bbb R^2, 0).$
From Theorem 5.1 the normal form of (III, I)$^1$ is
$(x+y, \gamma, f), f(u, 0)=u^{2},
\frac{\partial f}{\partial v}(0)\not=0.$
Assume that two normal forms of $(III, I)^{1}$ are topologically equivalent:
$$\CD
(\Bbb R,0) @<x+y << (\Bbb R^{2},0) @>\gamma>>
(\Bbb R^2,0) @>f>>  (\Bbb R, 0)\\
 @V h VV         @V H VV        @V K VV      @VV k V\\
(\Bbb R,0) @<X+Y << (\Bbb R^{2},0) @>\gamma>>
(\Bbb R^2,0) @>\widetilde{f}>> (\Bbb R, 0)\\
\endCD
\tag 5.2$$
where $h, H, K, k$ are homeomorphism germs.
We remark that 
if $f(u, v)=-v+u^2,$ then the parabola $f=0$ is the singular set
of the singular 3-web on $\{v>0\}$ associated with this normal form.
However clearly the type $(III, I)^{1}$ with $f(u, v)=-v+u^2$ is 
topologically equivalent to only oneself.
Hence we may suppose that $f(u, v)\not= -v+u^2.$ 
Then there exists a 3-web structure on $\{v>0\}.$
We remark that from (5.2)
$\frac{\partial f}{\partial v}(0)<0$ (resp. $>0$) if and only if
$\frac{\partial \widetilde{f}}{\partial v}(0)<0$ (resp. $>0$)
because the level curves $f^{-1} (0)$ with $\frac{\partial f}{\partial v}(0)
>0$ do not lie in the web domain $\{v>0\}.$ 
We see that 
$K|_{\{v>0\}}$ and $h$ are $C^{\infty}$-diffeomorphism germs
by the same argument as in (III, I)$^{0}.$

By Lemma 5.2 we have the following  $C^{\infty}$-commutative diagram
$$\CD
(\Bbb R,0) @<h(x+y) << (\Bbb R^{2},0) @>\gamma>>
(\Bbb R^2,0) @>f>>  (\Bbb R, 0)\\
 @|         @V \widehat{H} VV        @V \widehat{K} VV      @| \\
(\Bbb R,0) @<X+Y << (\Bbb R^{2},0) @>\gamma>>
(\Bbb R^2,0) @>f\circ\widehat{K}^{-1}>> (\Bbb R, 0)\\
\endCD
$$
where $\widehat{H}, \widehat{K}$ are $C^{\infty}$-diffeomorphism germs.
Therefore it is sufficient to show that $C^0$-
commutativity implies
$C^{\infty}$-commutativity in the following diagram:
$$\CD
(\Bbb R,0) @< x+y << (\Bbb R^{2},0) @>\gamma>>
(\Bbb R^2,0) @>f\circ\widehat{K}^{-1}>>  (\Bbb R, 0)\\
 @|         @V H\circ\widehat{H}^{-1} VV     @V K\circ\widehat{K}^{-1} VV      @VV k V \\
(\Bbb R,0) @<X+Y << (\Bbb R^{2},0) @>\gamma>>
(\Bbb R^2,0) @>\widetilde{f}>> (\Bbb R, 0) .\\
\endCD
\tag 5.3$$
We see that $H\circ\widehat{H}^{-1}=id$ on $(\Bbb R^2, 0)$
and  $K\circ\widehat{K}^{-1}=id$ on  $\{v\geq0\}.$

Now we consider the following two cases:

\noindent{\it Case $\,\frac{\partial f}{\partial v}(0)<0$} : 
By the same argument as the case (III, I)$^0 ,$
we see that $k$ is a $C^{\infty}$-diffeomorphism germ.
Define $C^{\infty}$-diffeomorphism germ $\widetilde{K} :
(\Bbb R^2, 0)\to (\Bbb R ^2, 0)$ by
$\widetilde{K}(u, v)=(u, \,\widetilde{f}(u, v))^{-1}\circ
(u, \, k \circ f\circ\widehat{K}^{-1}(u, v)).$
Then we have $\widetilde{K}=id$ on $\{v\geq 0\}$ and
$\widetilde{f}\circ\widetilde{K}
=k \circ f \circ\widehat{K}^{-1}.$
This completes the proof of rigidity 
in the case $\frac{\partial f}{\partial v}(0)<0$. 

\noindent
{\it Case  $\,\frac{\partial f}{\partial v}(0)>0$} : 
In this case $f^{-1}(0)$ is not in the 3-web domain $\{v>0\}.$
We see that $f$ restricted to $v$-axis ($v>0$) is a
$C^{\infty}$-diffeomorphism germ.
Then by the similar argument as in (III, I)$^0$ we see that 
$k|_{\{t>0\}}$ is a $C^{\infty}$-diffeomorphism germ.
Moreover from (5.3) we have 
$k(t)=\widetilde{f}(0, \varphi (t))$ for $t>0,$
where $\varphi (t)$ is the inverse function germ of
$t=f\circ \widehat{K}^{-1}(0, v).$
Let $\widetilde{k}(t)$ be a natural $C^{\infty}$-extension of $k|_{\{t>0\}}$
defined by $\widetilde{k}(t)=\widetilde{f}(0, \varphi(t))$
for all $t\in (\Bbb R , 0).$
Clearly $\widetilde{k}$ is a $C^{\infty}$-diffeomorphism germ.
Then similarly we can 
define $C^{\infty}$-diffeomorphism germ $\widetilde{K} :
(\Bbb R^2, 0)\to (\Bbb R ^2, 0)$ by
$\widetilde{K}(u, v)=(u, \,\widetilde{f}(u, v))^{-1}\circ
(u, \, \widetilde{k} \circ f\circ\widehat{K}^{-1}(u, v)).$
By the same argument as in the case $\frac{\partial f}{\partial v}(0)<0$  
we have the required result. 

This completes the proof of rigidity of case (III, I)$^1$. 
 \qed

\proclaim{5.3 Proof of rigidity for (VI, I)}
\endproclaim
In this case we suppose that $\gamma=(x, y^3 +xy).$
Denote by $\Delta$ the set 
$\{4u^3 +27v^2<0\}$
of inner points of the cusp. 
At first we prepare a crucial Proposition.
\proclaim{Proposition 5.3}
Let $(g_{i}, \gamma, f_i)$ be two normal forms of type $(VI, I)$
, where $g_{i}=y+\theta _i \circ\gamma, \theta_{i}\in \Cal M_2  (i=1, 2).$
Then $\theta_1 |_{\Delta} =\theta _2 |_{\Delta}$ and 
$f_1 |_{\Delta} = f_2 | _{\Delta}$ if and only if there exist $C^{\infty}$-
diffeomorphism germs $H, K : (\Bbb R ^2 , 0)\to (\Bbb R^2, 0)$
such that the following diagram commutes:
$$\CD
(\Bbb R,0) @<g_1 << (\Bbb R^{2},0) @>\gamma>>
(\Bbb R^2,0) @>f_1>>  (\Bbb R, 0)\\
 @|         @V H VV        @V K VV      @| \\
(\Bbb R,0) @<g_2 << (\Bbb R^{2},0) @>\gamma>>
(\Bbb R^2,0) @>f_2 >> (\Bbb R, 0).\\
\endCD
$$
\endproclaim
\demo{Proof}Necessity:
Assume that 
$\theta_1 |_{\Delta} =\theta _2 |_{\Delta}$ and 
$f_1 |_{\Delta} = f_2 | _{\Delta}.$ 
Define $C^{\infty}$-map germs $H_{i} : (\Bbb R^2, 0)\to (\Bbb R^2, 0)$
by $H_i=(f_{i}\circ\gamma, g_{i}), i=1,2.$
Since the level curve $f_{i}=0$ and the tangent cone
of the cusp are transversal, we see that $H_i$ is a $C^{\infty}$-diffeomorphism germ. 
Namely $H_i$ changes the 2-web $(f_i\circ\gamma, g_i)$ to the canonical grid.
Put $H=H_2 ^{-1}\circ H_1.$
Then we have $g_1 =g_2\circ H$ and $f_1\circ\gamma = f_2\circ\gamma
\circ H.$

On the other hand from the assumption clearly we have 
$H|_{\gamma^{-1}(\Delta)}=
id |_{\gamma^{-1}(\Delta)}.$
Then as was shown in [4]p.464 we see that $H$ is a $\gamma$-lowable, 
namely there exists a $C^{\infty}$-
diffeomorphism germ $K : (\Bbb R^2, 0)\to (\Bbb R^2, 0)$ such that 
$K\circ\gamma = \gamma\circ H.$
Hence  we see that these $H, K$ realize the required commutativity.

Sufficiency is a corollary of the argument of the below.
\qed
\enddemo

Assume that two diagrams of $(VI, I)$ are topologically equivalent:
$$\CD
(\Bbb R,0) @<g_1<< (\Bbb R^{2},0) @>\gamma>>
(\Bbb R^2,0) @>f_1>>  (\Bbb R, 0)\\
 @V h VV         @V H VV        @V K VV      @VV k V\\
(\Bbb R,0) @<g_2<< (\Bbb R^{2},0) @>\gamma>>
(\Bbb R^2,0) @>f_2>> (\Bbb R, 0)\\
\endCD
$$
where $g_i =y_i + \theta _i \circ \gamma (i=1, 2)$ and
 $h, H, K, k$ are homeomorphism germs.

\proclaim{Step 1} We may suppose that $h=id.$
\endproclaim
There is a 4-web structure in $\Delta$ associated with $(VI, I).$
Hence by rigidity of webs $K|_{\Delta}$ is a $C^{\infty}$-diffeomorphism germ.
Then as was shown in [4]p.469
(which is the similar argument as in the previous cases) 
we see that $h$ is a $C^{\infty}$-diffeomorphism germ.
Hence $(h\circ g_1, \gamma, f_1)$ is $C^{\infty}$-equivalent to
$(g_1, \gamma, f_1).$
On the other hand by [4, Proposition 1 (p. 462)] we have the following
$C^{\infty}$-commutative diagram
$$\CD
(\Bbb R,0) @<h\circ g_1 << (\Bbb R^{2},0) @>\gamma>>
(\Bbb R^2,0) @>f_1 >>  (\Bbb R, 0)\\
 @|         @V \widehat H VV        @V \widehat K VV      @| \\
(\Bbb R,0) @< y+\theta\circ\gamma << (\Bbb R^{2},0) @>\gamma>>
(\Bbb R^2,0) @>\widehat{f}=f_{1}\circ\widehat{K} ^{-1}>> (\Bbb R, 0)\\
\endCD
$$
for some $\theta\in\Cal M_2$,
where $\widehat {H}, \widehat{K}$ are $C^{\infty}$-diffeomorphism  germs.
Namely we have the following  $C^{0}$-commutative diagram:
$$\CD
(\Bbb R,0) @<y+\theta\circ\gamma << (\Bbb R^{2},0) @>\gamma>>
(\Bbb R^2,0) @>\widehat{f} >>  (\Bbb R, 0)\\
 @|       @V H\circ\widehat{H}^{-1} VV        @V K\circ\widehat{K}^{-1} VV      @VV k V\\
(\Bbb R,0) @<g_2<< (\Bbb R^{2},0) @>\gamma>>
(\Bbb R^2,0) @>f_2>> (\Bbb R, 0).\\
\endCD
$$
Since the level curve $\widehat{f}=0$ and the tangent cone
of the cusp are also transversal, we can suppose that $\widehat{f}(u, 0)=u$ by the
coordinate change,  
which is the inverse function of $\widehat{f}(u, 0),$
in the target of $\widehat{f}.$
Hence it is sufficient to show the case $h=id.$

\proclaim{Step 2} $\theta _1 |_{\Delta}=\theta_{2}|_{\Delta},\,
f_{1}|_{\Delta}=f_{2}|_{\Delta}$ and apply Proposition 5.3
\endproclaim

By the same argument as in [4]p.470, due to $h=id$
we see that $H|_{\gamma^{-1}(\Delta)}=id|_{\gamma^{-1}(\Delta)}$,
$K|_{\Delta}=id|_{\Delta}.$
Hence $\theta_1 =\theta_2$ on $\Delta.$
Moreover we have $k(s)=s$ for $s<0.$
Hence $f_1 =f_2$ on $\Delta.$ 
Then applying the necessity part of Proposition 5.3,
the proof of rigidity for $(VI, I)$ is completed.
\qed

\head{Acknowledgments}
\endhead
\eightpoint
\baselineskip=12pt

\noindent 
A part of this work was done when the second author visited GETODIM
Institut de Math\'{e}matiques Universit\'{e} de Montpellier II,
supported by Research Fellowships of Japan Society for the Promotion
of Science, and this work was completed when he is in ICMSC-Instituto de
Ci\^{e}ncias Matem\'{a}ticas de S\~{a}o Carlos, Universidade de 
S\~{a}o Paulo, supported by FAPESP grants.
He would like to thank GETODIM and ICMSC for their hospitalities.
Also the second author would like to thank
Prof.~S.~Izumiya for his valuable advices and constant encouragements,
Prof.~S.~Izumi for his useful discussion
and Prof.~M.~A.~S.~Ruas for her helpful comments.

\tenpoint

\midspace{0.5cm}
\baselineskip=12pt

\noindent
Jean Paul Dufour

\noindent
Getodim Institut de Math\'{e}matiques, Universit\'{e} de Montpellier II

\noindent
Pl. E. Bataillon 34060 Montpellier Cedex, France

\noindent
E-mail : dufourj{\char'100}darboux.math.univ-montp2.fr

\bigpagebreak

\noindent
Yasuhiro Kurokawa

\noindent
Instituto de Ci\^{e}ncias Matem\'{a}ticas de S\~{a}o Carlos

\noindent
Universidade de S\~{a}o Paulo, Caixa Postal 668, CEP 13560-970, S\~{a}o Carlos, SP, Brazil

\noindent
E-mail : kurokawa{\char'100}icmc.sc.usp.br

\end